\documentclass[12pt]{amsart}
\usepackage{amsmath,amssymb,amsfonts,amsthm,xspace}
\usepackage{graphicx}

\newtheorem{theorem}{Theorem}

\newtheorem{lemma}[theorem]{Lemma}
\newtheorem{corollary}[theorem]{Corollary}

\newcommand{\be}{\begin{equation}}
\newcommand{\ee}{\end{equation}}
\newcommand{\Z}{{\mathbb Z}}
\newcommand{\R}{{\mathbb R}}

\newcommand{\M}{{\mathcal M}}
\newcommand{\C}{{\mathbb C}}

\newcommand{\G}{{\mathcal G}}
\newcommand{\SL}{\mathrm{SL}_2({\mathbb C})}
\newcommand{\SU}{\mathrm{SU}_2}
\newcommand{\GL}{\mathrm{GL}_2({\mathbb C})} 
\newcommand{\Tr}{\text{Tr}}
\newcommand{\Qdet}{\mathrm{Qdet}}
\newcommand{\old}[1]{}
\renewcommand{\b}{\mathbf{b}}
\newcommand{\w}{\mathbf{w}}
\newcommand{\f}{\mathbf{f}}
\newcommand{\Pf}{\mathrm{Pf}}
\newcommand{\V}{{\mathcal V}}
\newcommand{\B}{{\mathcal B}}
\newcommand{\W}{{\mathcal W}}
\renewcommand{\H}{{\mathbb H}}
\newcommand{\eps}{{\epsilon}}
\renewcommand{\Im}{{\text{Im}}}
\renewcommand{\Re}{{\text{Re}}}

\title{Conformal invariance of loops in the double-dimer model}
\author{Richard Kenyon }

\begin{document}

\begin{abstract} 
The dimer model is the study of random dimer covers (perfect matchings) of a graph.
A double-dimer configuration on a graph $\G$ is a union of two dimer
covers of $\G$. We introduce
quaternion weights in the dimer model
and show how they can be used to study the
homotopy classes (relative to
a fixed set of faces) of loops in the double dimer model on a planar graph. As an application we
prove that, in the scaling limit of the ``uniform" double-dimer model on $\Z^2$
(or on any other bipartite planar graph conformally approximating $\C$),
the loops are conformally invariant. 

As other applications we compute the exact
distribution of the number of topologically nontrivial
loops in the double-dimer model on a cylinder
and the expected number of loops surrounding two faces of a planar graph.
\end{abstract}
 
\maketitle

\section{Introduction}

A {\bf dimer cover}, or perfect matching, of a graph $\G$ is a set of edges
with the property that each vertex
is the endpoint of exactly one edge. 
A {\bf double-dimer configuration} is a union of two dimer covers, or
equivalently a set of even-length simple loops and doubled edges with the 
property that every vertex is the endpoint of exactly two edges (which may be doubled).

The uniform dimer model (the uniform measure on dimer covers)
on the square grid $\Z^2$ has been the subject of much research,
starting with the exact enumeration of \cite{Kast, TF} and culminating in
the limit shape theorems in, successively,
\cite{CKP, KOS, KO}, and conformal invariance results in \cite{K.gff}.
It was conjectured in \cite{KW} (see also \cite{Sheff}) that the loops in the double-dimer model also
have conformally-invariant scaling limits, and are related to $SLE_4$.
We prove here the conformal invariance of the loops.

If a graph $\G$ has positive edge weights $\nu\colon E\to\R_{>0}$, 
there is a natural probability measure $\mu_{\nu}$ on dimer covers
which gives a dimer cover a probability proportional to the product of its edge weights.
From $\mu_\nu\otimes\mu_\nu$ we similarly get probability measures on double-dimer configurations.
More general double-dimer probability measures can be constructed by taking 
$\mu_{\nu_1}\otimes\mu_{\nu_2}$ for different $\nu_1,\nu_2$, or even
$\nu_2=\bar\nu_1$ when $\nu_1$ is complex-valued (in this case $\mu_{\nu_1}\otimes\mu_{\nu_2}$ may be a probability measure on double-dimer configurations,
thought of as collections of loops,
even though $\mu_{\nu_1}$ and $\mu_{\nu_2}$ are not probability measures). 

The main innovation in this paper is the introduction of quaternionic (or more generally
$\SL$) weights
for the double-dimer model, a tool which is not available for the single-dimer model.
The Kasteleyn theorem on evaluating the weighted sum of dimer covers via 
determinants can
be extended to apply in this case; one must use the so-called quaternion-determinant,
or $q$-determinant
\cite{Dyson, Mehta}
of a self-dual quaternionic matrix. The $q$-determinant of the quaternionic Kasteleyn matrix is then the partition function of the double-dimer model 
(Theorem \ref{main1} below).

Another important tool in the proof of conformal invariance is a theorem of Fock and Goncharov \cite{FG}, Theorem \ref{FGthm} below,
connecting the set of simple closed curves (or more generally laminations)
on surfaces with the $\SL$ representation theory of its fundamental group.

On a multiply-connected planar domain $U$ we show that, for a graph 
conformally approximating
$U$ (in the sense that discrete harmonic functions converge to 
harmonic functions on $U$), the distribution of the homotopy types of the
homotopically nontrivial double dimer loops only depends on the 
conformal type of $U$. 
See Corollary \ref{mainci} for the case when the graph is
essentially $U\cap \eps\Z^2$, and Section \ref{othergraphs} for the extension to other
graphs.

As a limiting case of the above we can take $U$ to be a bounded simply connected region from which a finite number of points  
$z_1,\dots,z_m$ have been removed; in this case we show that the distribution
of the homotopy types of the nonperipheral, homotopically nontrivial loops in the double dimer model
only depends on the conformal type of $U$. (A loop is peripheral if it is isotopic to a small loop surrounding
exactly one of the punctures $z_i$.) See Section \ref{peripheral}.

We give two other applications of these ideas: a computation of
the exact distribution of the number of 
topologically non-trivial loops in a double-dimer configuration
on the square-grid on a cylinder (see Section \ref{annulus}),
and the expected number of double-dimer loops surrounding simultaneously two
points $z_1,z_2$ (see Section \ref{twofaces}).

In our main proof we approximate $U$ with a 
discrete graph $\G_\eps\subset\eps\Z^2$ with rather special boundary conditions called 
Temperleyan boundary conditions \cite{K.ci},
(See Figure \ref{Geps}). This issue of choice of boundary conditions is a delicate one.
In fact, for naive choices of boundary conditions the conformal invariance will not hold,
see \cite{Kenyon.honeycomb}. There are nonetheless
many other ``good" choices of boundary
conditions and it is possible to show with the techniques in this paper 
that the conformal invariance holds for these
graphs as well. Our choice in the current paper was designed to make the proof
as simple as possible. The technical difficulty in the general
case is showing convergence of $K^{-1}$, the inverse Kasteleyn operator.
In the current case we can write $K^{-1}$ in terms of the discrete 
Dirichlet Green's function whose convergence properties are well known.
For other boundary conditions the corresponding discrete Green's functions are
(usually) less well understood.

\section{Vector bundles and connections}

Let $\G=(V,E)$ be a graph and $W$ a vector space.
A {\bf $W$-bundle} on $\G$ is a vector space $W_v$ isomorphic to $W$ associated
to each vertex $v$ of $\G$. The {\bf total space} of the bundle is 
$W_{\G}= \oplus_{v\in V} W_v$; a {\bf section} is an element of $W_{\G}$.
A {\bf connection} $\Phi=\{\phi_e\}_{e\in E}$ on a vector bundle
is the choice of an isomorphism
$\phi_{vv'}\colon W_v\to W_{v'}$ for every edge $e=vv'$, 
with the property
that $\phi_{-e}=\phi_{e}^{-1}$, where $-e$ is the edge $e$ with the reverse orientation.
The map $\phi_{vv'}$ is the {\bf parallel transport}
of vectors at $v$ to vectors at $v'$.

Two connections $\Phi$ and $\Phi'$ are {\bf gauge equivalent} if there
are isomorphisms $\psi_v\colon W_v\to W_v$ such that 
$\psi_{v_2}\phi'_{v_1v_2}=\phi_{v_1v_2}\psi_{v_1}$ for all adjacent $v_1,v_2$.

If $\gamma$ is an oriented closed path on $\G$ starting at $v$, 
the {\bf monodromy} of the connection
along $\gamma$ is the product of the parallel transports along $\gamma$;
it is an element of $\text{End}(W_v)$. Gauge equivalent connections give
conjugate monodromies; changing the starting point along $\gamma$
also conjugates the monodromy. 

If $\G$ is embedded on a surface (in such a way that faces are contractible),
a {\bf flat connection} is a connection with trivial monodromy around faces,
and thus around any contractible curve. If $\Phi$ is a flat connection, 
given any contractible union of faces $S$, it is not hard to see that
one can choose a connection gauge equivalent
to $\Phi$ which is locally trivial on $S$, that is, is the identity on each edge in $S$. 

It is convenient to represent a flat connection as in Figure \ref{Geps}: take a set of simple 
closed paths in the dual graph which support cocycles generating the cohomology of the surface
(in other words, cutting the edges crossing these paths results in a contractible surface)
and take flat connection supported on the edges crossing these paths. Any flat connection is equivalent
to one of this form. We call such a dual path a zipper. 

In this paper we use only $1$- or $2$-dimensional complex vector bundles,
with parallel transports in $\C^*$ or $\SL$. In fact for the purposes of 
studying the double-dimer measure it suffices
to take unitary connections, with structure group ${\text U}_1$ or $\SU$.

\section{Dimer model}

\subsection{The single-dimer model}
\subsubsection{Definition}
Let $\G$ be a graph and $\nu:E\to\R_{>0}$ a positive real weight function on the edges.
Let $\M=\M(\G)$ be the set of dimer covers (perfect matchings) of $\G$.
For $m\in\M$ define $\nu(m)=\prod_{e\in m}\nu(e)$ to be its weight. 
We define a probability measure $\mu$ on $\M$ where the probability of $m$ is proportional to $\nu(m)$. The constant of proportionality is $1/Z$ where
$Z=\sum_{m\in\M}\nu(m)$ is
called the {\bf partition function}.

Note that if $\nu'$ is a different weight function,
obtained from $\nu$ by multiplying the edge
weights at a given vertex by a nonzero constant $\lambda$, then $\mu$ is unchanged,
since every dimer cover has weight multiplied by $\lambda$.
Compositions of such operations are called {\bf gauge transformations}. Two
weight functions are {\bf gauge equivalent} if they are related by a gauge transformation.

\subsubsection{Line bundle interpretation}

Now suppose $\G$ is bipartite, with black vertices $\B$ and white vertices $\W$.
Define $\Phi=\Phi_\nu$, a connection on a line bundle on $\G$, by taking $W=\C^1$
and $\phi_{e}=\nu(e)$
if $e$ is directed from black to white, that is,
$\phi_e$ is multiplication by $\nu(e)$ (we must then have 
$\phi_{-e}=\nu(e)^{-1}$).
If $\nu_1$ and $\nu_2$ are gauge equivalent as weight functions, then their associated
connections $\Phi_{\nu_1},\Phi_{\nu_2}$ are gauge equivalent as bundles.
Thus $\mu$ depends only on the gauge equivalence class of $\Phi$; 
gauge transformations do not change $\mu$.

\subsubsection{Kasteleyn matrix}

If $\G$ is planar and bipartite, Kasteleyn showed \cite{Kast} that the partition 
function $Z$ is
the determinant of the {\bf Kasteleyn matrix}
$K=(K(w,b))_{w\in\W,b\in\B}$ whose rows index white vertices and columns index black vertices, and $K(w,b)=\pm\nu_{wb}=\pm \phi_{wb}$. Here the sign is chosen according
to the ``Kasteleyn sign condition" for bipartite graphs, that is, so that
the number of minus signs around a face of length $\ell$ is $\frac{\ell}2+1 \bmod 2$. 
We have $Z=|\det K|$. 

In the case $\G$ is planar but not bipartite,
there is a generalization of this result also due to Kasteleyn: $Z$ is the Pfaffian of an antisymmetric
matrix $K=(K(v,v'))_{v,v'\in V}$ (also called Kasteleyn matrix), indexed by all the vertices of $\G$, with $K(v,v')=\pm\nu(vv')$. Here the sign is chosen by orienting
the edges of $\G$ in such a way that each face has an odd number of clockwise-oriented
edges.

\subsection{Double dimer model}

\subsubsection{Definition}
Let $\G$ be a general (not necessarily bipartite) graph
and $\nu_1,\nu_2:E\to\C$ two weight functions, not necessarily real. 
Given a pair $(m_1,m_2)\in\M^2$, we associate a weight
$\nu(m_1,m_2)=\nu_1(m_1)\nu_2(m_2)$ where $\nu_i(m)$ are defined as above.
The partition function $Z_{dd}$ is defined as
\be\label{Z1Z2}
Z_{dd}= Z(\nu_1)Z(\nu_2)=\sum_{(m_1,m_2)\in \M^2}\nu(m_1,m_2) .\ee
 
It is natural to group configurations according to the set of loops they form.
Let $\Omega=\Omega(\G)$ be the set coverings of $\G$ by 
collections of edges which form cycles of even length
or doubled edges (so that each vertex is the endpoint of exactly two edges,
which may be the same edge). We call $\Omega$ the set of
{\bf double dimer configurations} on $\G$. 
To $\omega\in\Omega$ we associate a weight 
$$\nu(\omega)=\prod_{\text{doubled edges}}\nu_1(e)\nu_2(e)\prod_{\text{cycles}}(w_1+w_2)$$
where for a cycle $\gamma=(v_1\to v_2\to\dots\to v_{2n}\to v_1)$ we associate two weights 
$$w_1(\gamma)=\nu_1(v_1v_2)\nu_2(v_2v_3)\nu_1(v_3v_4)\nu_2(v_4v_5)\dots$$
and
$$w_2(\gamma)=\nu_2(v_1v_2)\nu_1(v_2v_3)\nu_2(v_3v_4)\nu_1(v_4v_5)\dots.$$

\begin{lemma}
$Z_{dd}=Z(\nu_1)Z(\nu_2)=\sum_{\omega\in\Omega}\nu(\omega).$
\end{lemma}

\begin{proof}
A single element $\omega\in\Omega$ with $k$ cycles corresponds to $2^k$
pairs $(m_1,m_2)\in\M^2$: each cycle in $\omega$ can be partitioned in two ways
to get two dimer covers of it. 
\end{proof}

If we are interested in constructing a probability measure on $\Omega$, 
we must (somehow) choose $\nu_1,\nu_2$ so that the 
weights of configurations are real and nonnegative. The reality
can be guaranteed if we take for example $\nu_2=\bar\nu_1$, since then
$w_2=\bar w_1$ for any cycle. Positivity can then 
be achieved if the arguments of the weights
are sufficiently small. 

\subsubsection{Computation}
If $\G$ is planar, from (\ref{Z1Z2}) we have 
\be\label{Z1} Z_{dd}=\Pf K_1\Pf K_2\ee where $K_1,K_2$ are antisymmetric Kasteleyn matrices associated to $\nu_1,\nu_2$ respectively. 

If $\G$ is bipartite and planar,
we can instead use a single matrix $K=(K(v,v'))_{v,v'\in V}$ indexed by all vertices,
where
$K(w,b)=\pm\nu_1(b,w)$ and $K(b,w)=\pm\nu_2(w,b)$. The signs are again determined
by the Kasteleyn condition for bipartite graphs mentioned above. In this case we have
\be\label{Z2}
Z_{dd}=\det K.
\ee
Equations (\ref{Z1}) and (\ref{Z2}) are special cases (for antidiagonal matrices and diagonal matrices, respectively) of the construction below.

\subsection{$\SL$-bundle case}

Let $\G$ be a bipartite graph and $\nu\colon E\to\R_{>0}$ a positive real weight function
on the edges.
Let $\Phi$ be an $\SL$-connection on a $\C^2$-bundle on $\G$, that is, a connection
where the parallel transports are in $\SL$.

The weight of $\omega\in\Omega$ is now defined to be
\be\label{connectionweight}
\nu(\omega) = \prod_{e} \nu(e)\prod_{\text{cycles}}\Tr(w),
\ee
where the first product is over edges of $\omega$ 
(and doubled edges are counted twice)
and the second product is over nontrivial cycles (that is, not the doubled edges), 
$w$ being the monodromy of the cycle
in one direction or the other (note that $\Tr(w)=\Tr(w^{-1})$).  

In order to associate a probability measure it is necessary that the trace
of the monodromy on any loop which occurs in some $\omega\in\Omega$ be nonnegative. 
This can be arranged for example if the parallel transports $\phi_{vv'}$ 
are all sufficiently close to the identity.

Associated to an $\SL$ connection on a bipartite planar graph is a Kasteleyn
matrix $K$ with entries in $\GL$, defined as follows. We define 
$K(v,v')=K_\nu(v,v')\phi_{vv'}$ where $K_\nu$ is the usual weighted signed
Kasteleyn matrix, and $\phi_{vv'}$ are the parallel transports written in the standard
basis for $\C^2$.

Note that the matrix $K$ is 
{\bf self-dual}, that is $K(v,v')=K(v',v)^*$
where by ${}^*$ we mean the ``$q$-conjugate": if $A=\left(\begin{matrix}a&b\\c&d\end{matrix}\right)$ then $
A^* = \left(\begin{matrix}d&-b\\-c&a\end{matrix}\right).$

\begin{theorem}\label{main1}
We have $Z_{dd}=\sum_{\omega\in\Omega}\nu(\omega)= \Qdet K$.
\end{theorem}

Here $\Qdet$ is the quaternion determinant of a self-dual matrix, 
defined by:
\be\label{qdetdef}
\Qdet K = \sum_{\sigma\in S_n}\text{sgn}(\sigma)\prod_{\text{cycles }C}\frac12\Tr(K_C),
\ee
where the product is over cycles $C$ of $\sigma$, and $K_C$ is the product of the elements of $K$ along the cycle $C$, that is, if $C$ is the cycle 
$$v_1\to v_2\to\dots v_n\to v_1$$ then $K_C=K_{v_1v_2}K_{v_2v_3}\dots K_{v_nv_1}$.

\begin{proof}
Each summand in (\ref{qdetdef}) for $\Qdet(K)$ is zero unless $\sigma$ maps nearest neighbors to 
nearest neighbors, and thus is a double-dimer configuration; because $K$ is bipartite
all cycles must have even length. 
The weight of a nonzero term is equal to the product of its edge weights
(with doubled edges counting twice) times the product of $1/2$ the trace
of the monodromy on all nontrivial (length bigger than $2$) cycles. When we account for 
both orientations of each nontrivial cycle, its contribution is the trace of its monodromy.

It remains to verify that each nonzero term has the same sign. This is true
by Kasteleyn's theorem (equation (\ref{Z2}) above)
for the trivial bundle, and thus (by for instance
deforming continuously to the case of a nontrivial bundle) remains
true in general.
\end{proof}

If $K$ is an $n\times n$ self-dual matrix with entries in $\GL$, let $\tilde K$ be the 
$2n\times 2n$ matrix with entries in $\C$, obtained by replacing each entry by its $2\times 2$ block.
Mehta \cite{Mehta2}, see also \cite{Dyson},
proved that $\Qdet K$ is equal to the Pfaffian of the (antisymmetric) matrix 
${\mathcal Z}\tilde K$, 
where ${\mathcal Z}$ is the $2n\times 2n$
matrix with blocks $\left(\begin{matrix}0&1\\-1&0\end{matrix}\right)$ on the diagonal
and zeros elsewhere. 

Since $\G$ is bipartite, listing white vertices first, $K$ has the form 
$K=\left(\begin{matrix}0&M\\M^*&0\end{matrix}\right)$ where $M$ is the matrix
with $\GL$-entries with
rows indexing white vertices and columns indexing black vertices. Thus
$\Qdet K = \det\tilde M$ where $\tilde M$ is the 
$n\times n$ matrix obtained by replacing each entry in $M$ by its $2\times 2$ block.

Later we will use the fact that the inverse of a self-dual matrix $K$
of nonzero $\Qdet$ is well-defined, is both a left- and right-inverse, 
and is self-dual, see \cite{Dyson}.

\label{diagonalmtx}
Note that if we choose $\phi_e= \left(\begin{matrix}\lambda&0\\0&\lambda^{-1}\end{matrix}\right)$ where $\lambda=\sqrt{\frac{\nu_1(e)}{\nu_2(e)}}$
and $\nu(e)=\sqrt{\nu_1(e)\nu_2(e)}$ then (\ref{Z2}) is a special case of Theorem \ref{main1}. 

The construction in this section also works for nonbipartite graphs if we can
arrange that the trace of the monodromy around every odd-length loop is $0$
(then configurations with odd-length loops have zero weight).
For example if we choose $\phi_e= \left(\begin{matrix}0&\lambda\\-\lambda^{-1}&0\end{matrix}\right)$ where $\lambda=\sqrt{\frac{\nu_1(e)}{\nu_2(e)}}$
and $\nu(e)=\sqrt{\nu_1(e)\nu_2(e)}$ then (\ref{Z1}) is a special case of Theorem 
\ref{main1} for nonbipartite graphs. 
\medskip

\noindent{\bf Question:} Are there other, more interesting, connections
in which the trace of the odd-length-loop monodromies are zero?

\section{Topology of paths}

An {\bf noncontractible simple closed curve} on a surface $\Sigma$ is a simple closed curve
which does not bound a disk.
A {\bf finite lamination} on a surface $\Sigma$, possibly with boundary,
is a finite set $L=\{\gamma_1,\dots,\gamma_k\}$
of pairwise disjoint noncontractible simple closed 
curves $\gamma_i$, considered up to isotopy.
By this we mean that we consider two finite laminations to be equivalent if they are isotopic,
that is, they have the same number of loops in the same homotopy 
classes.

Let $(\G,\nu,\Sigma)$ be a triple consisting of a
bipartite graph $\G$, with edge weights $\nu$, embedded on a surface $\Sigma$. 
The loops in a double-dimer cover $\omega\in\Omega(\G)$ may or may not be contractible. 
We associate to $\omega$ the finite lamination $\lambda(\omega)$ 
which is the (possibly empty)
set of noncontractible cycles, that is, those cycles which do not bound
disks in $\Sigma$. 

Let $\mu_0$ be the probability measure  
on double-dimer configurations $\Omega$ defined for $\omega\in\Omega$ by
$$\mu_0(\omega)=\frac{2^k}{Z}\prod_{e\in\omega} \nu(e),$$
where $Z$ is a normalizing constant,
$k$ is the number of nontrivial cycles (contractible or not)
and again doubled edges are counted twice in the product (but do not count as cycles).
This is the weight (\ref{connectionweight})
associated with the trivial connection. We refer to $\mu_0$ as the {\bf natural} measure on double dimer configurations (since it arise from the uniform measure on single-dimer configurations). 
We wish to study the 
$\mu_0$-probability of any finite lamination on $\Sigma$, that is, the probability
that a $\mu_0$-random configuration $\omega$
has $\lambda(\omega)$ of a specified type.

Recall that a flat connection on $(\G,\Sigma)$ is
a connection whose monodromy around every contractible 
loop is the identity. The set of flat $\SL$-connections modulo
gauge transformations has the structure of an algebraic
variety $X=X(\Sigma)$, isomorphic to the {\bf representation variety} of homomorphisms of 
$\pi_1(\Sigma)$ into $\SL$, modulo conjugation.

Let $\Phi=\{\phi_e\}$ be a flat connection.
By Theorem \ref{main1}, the partition function $Z_{dd}=Z_{dd}(\Phi)$,
which is a polynomial function of the entries of the $\phi_e$,
is a nonnegative linear combination of terms $\prod_{\gamma\in L} \Tr(w_\gamma)$ where
$L$ runs through all possible finite laminations of $(\G,\Sigma)$. In particular
it is a function on $X$.

Let $\V=\V(\G,\Sigma)$ be the vector space of functions of 
$\Phi=\{\phi_e\}$ which are polynomial
in the matrix entries and invariant under gauge transformations.
This is a vector space of infinite dimension but finite in each degree.
Fock and Goncharov proved 
\begin{theorem}[\cite{FG}]\label{FGthm} Suppose $\Sigma$ has nonempty boundary.
The elements $\prod_{\gamma\in L} \Tr(w_\gamma)$, as elements of $\V$, are linearly independent,
and in fact form a basis for $\V$ as $L$ runs over all finite laminations.
Moreover there is a natural inner product on $\V$, coming from Haar measure
on $\SL$, under which the functions $\prod_{\gamma\in L} \Tr(w_\gamma)$ are orthogonal.
\end{theorem}

In particular one can extract, via an integral over the representation variety $X(\Sigma)$,
the coefficient in $Z_{dd}$ of any given finite lamination. When the surface has nonempty
boundary, $\pi_1(\Sigma)$ is a free group and so 
$X(\Sigma)$ is a just a product of copies of $\SL$. The corresponding integrals
are integrals against spherical harmonics times the Haar measure on $(\SU)^n$.

\section{Scaling limit of the double-dimer path}

\subsection{Scaling limit}\label{setup}

Let $U\subset\C$ be a bounded multiply connected
domain with boundary consisting of a finite number of piecewise smooth 
curves $C_0,\dots,C_m$, with $C_0$ being the outer curve. 
Let $z_i\in C_i$ be fixed points, one on each boundary component.
We assume (to simplify the proof of convergence of the Green's function)
that the boundary of $C_i$ is flat and horizontal in a neighborhood of each $z_i$.
Let $\gamma_i$ be a simple path in the interior of $U$ from $z_0$ to $z_i$. 
We suppose these $\gamma_i$ are pairwise disjoint except at $z_0$.

\begin{figure}[htbp]
\includegraphics[height=2.in]{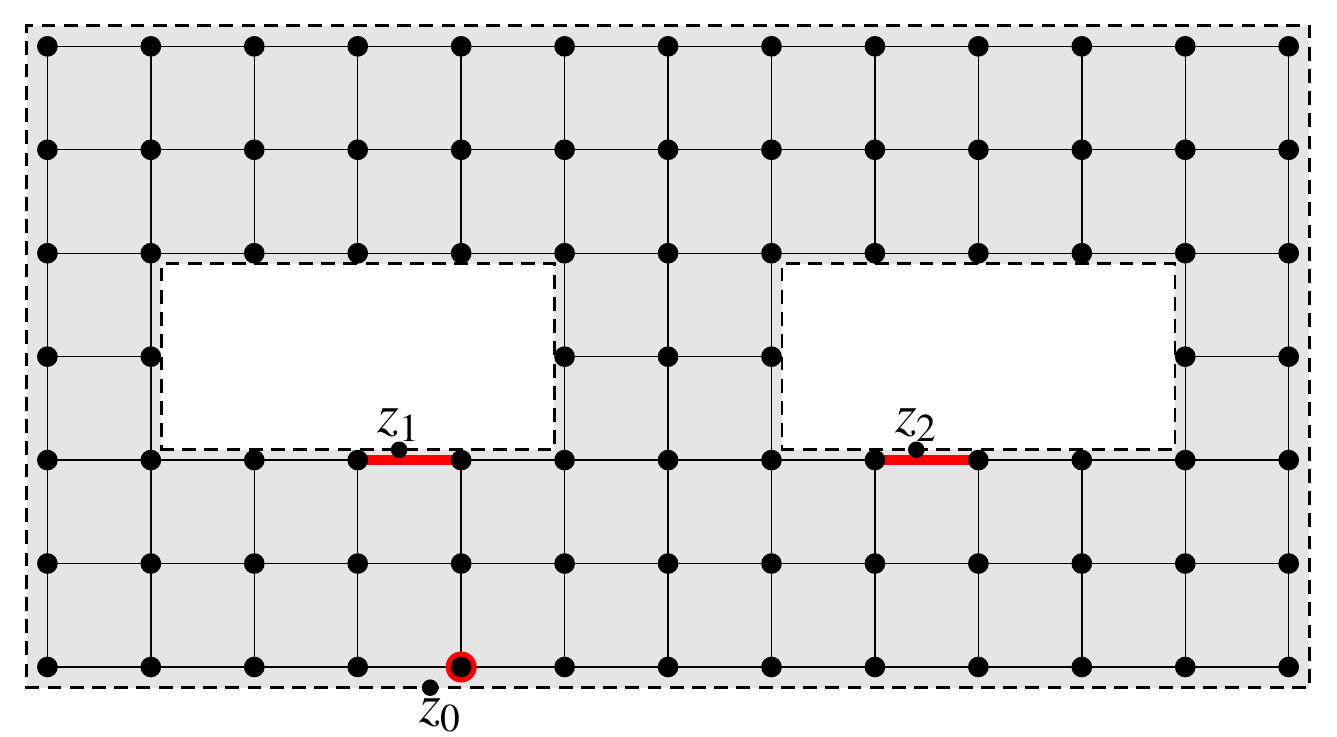}
\caption{\label{Ueps}The region $U$ and 
graph $U_\eps$ with marked edges $e_i$ and vertex $x_0$.}
\end{figure}

For $\eps>0$ let $U_\eps$ be the graph with vertices $U\cap\eps\Z^2$
and edges connecting points at distance $\eps$. For each $j>0$ let $e_j$ be an edge 
of $U_\eps$ near $z_j$ (that is, converging to $z_j$ as $\eps\to0$),
whose vertices lie on the face containing $C_j$. Let $x_0$ be
a vertex near $z_0$. See Figure \ref{Ueps} for an example.

\begin{figure}[htbp]
\includegraphics[height=2.in]{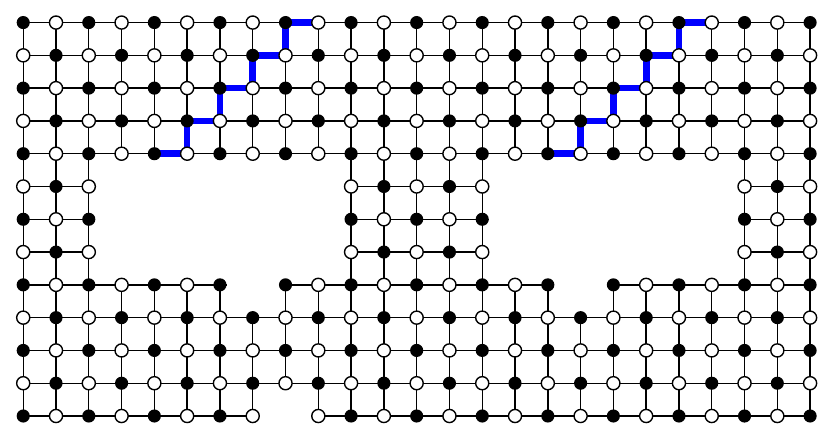}
\caption{\label{Geps}The associated graph $G_\eps$. The ``zippers" are the thicker edges 
(in blue).}
\end{figure}

Let $\G_\eps$ be a bipartite graph with a black vertex for each vertex (except $x_0$)
and square face
of $U_\eps$ (faces do not include the complementary components), and a white 
vertex for every edge except for
the edges $e_i$. Edges in $\G_\eps$ connect nearest neighbors, that is,
``edge" vertices to ``vertex" vertices and 
``edge" vertices to ``face"
vertices if the corresponding elements are adjacent. See Figure \ref{Geps}.
The graph $\G_\eps$ is bipartite and has dimer covers (these are in fact in bijection
with spanning trees of $U_\eps$ rooted at $x_0$ 
with a certain property: their dual trees 
rooted on the outer face have branches from the dual vertices in the $C_i$ which
lead in the direction of the $e_i$;  see \cite{KPW}. 

Let $\mu_\eps$ be the natural measure on double-dimer configurations 
$\Omega(\G_\eps)$ of $\G_\eps$,
that is, in which a configuration with $c$ loops 
has a probability proportional to $2^c$.

For a double-dimer configuration $\omega\in\Omega(\G_\eps)$ 
let $L=L(\omega)$ be the finite lamination of $U$ consisting of the noncontractible
loops of $\omega$.

We prove that the $\mu_\eps$-distribution of $L$ 
converges as $\eps\to 0$ and the limit only depends on the conformal type of 
the marked surface $(U,\{z_0,\dots,z_m\})$.

\subsection{Conformal invariance}

Let $\Phi$ be a flat $\SL$-connection on $(\G_\eps,U)$.
Let $A_i\in\SL$ be the monodromy of a path in $\G_\eps$
isotopic to one traversing $\gamma_i$, then running
counterclockwise
around the boundary of $C_i$, and then traversing $\gamma_i$ back towards $z_0$.

Let $Z_{dd}^{(\eps)}=Z_{dd}^{(\eps)}(\Phi)$ be the double-dimer partition
function; for each $\eps$ it is a polynomial function of the entries of the $A_i$.
Let 
\be\label{Flimit}
F(\Phi) = \lim_{\eps\to0}\frac{Z_{dd}^{(\eps)}(\Phi)}{Z_{dd}^{(\eps)}(\text{Id})}
\ee
where $\text{Id}$ is the trivial connection.

\begin{theorem}\label{CI}
The above limit exists and the function $F(\Phi)$ is conformally invariant,
that is, depends only on the conformal type of $(U,\{z_1,\dots,z_m\})$.
\end{theorem}

The proof will show that $F(\Phi)$ is an integral of certain products of Green's functions
on $U$. The geometric meaning of $F$ is not clear in general; 
however for the case when $U$ is an annulus see section \ref{annulus} below. 

Combined with Theorem \ref{FGthm} we have the following corollary. 
\begin{corollary}\label{mainci} The limiting distribution of the 
lamination $L$ is conformally invariant.
\end{corollary}

\noindent{\bf Proof of Theorem \ref{CI}.}
We can realize $\Phi$ by choosing $\phi_{vv'}$ to be the identity on all edges
except for the edges crossing a set of disjoint simple paths from the $C_i$ to the 
outer boundary $C_0$,
as in Figure \ref{Geps} (in blue). For simplicity we choose these edges
so that, oriented from the boundary to $f_i$, each crossing edge has a white 
vertex on the left and black vertex on the right. We call these collections of edges zippers.
Let $E_i$ be the set of edges in the $i$th zipper.

Let $\{\Phi_t\}_{t\in[0,1]}$ be a smooth one-parameter family of flat connections
supported on the zippers with $\Phi_0$ being the identity and $\Phi_1=\Phi$.
It suffices to show that $\frac{d}{dt}\log Z_{dd}^{(\eps)}(\Phi_t)$ is a conformally invariant
quantity plus an error tending to zero as $\eps\to0$.

Fix $\eps>0$. Let $K_t$ be a Kasteleyn matrix associated to $\Phi_t$ as above, with 
a fixed sign convention independent of $t$. 
Let $\tilde K_t$ be the matrix obtained from $K_t$ by replacing
each $\SL$ entry with its $2\times2$ block of complex numbers; thus
$\tilde K_t$ is a matrix of twice the dimension; 
to each vertex $v\in\G_\eps$ there are two rows and two columns
$v^{(1)}$ and $v^{(2)}$ of $\tilde K_t$. Under an infinitesimal change $t\to t+\delta$, 
$\tilde K_{t+\delta}=\tilde K_t+\delta \tilde S_t$ for some
matrix $\tilde S_t$ supported on the zippers.
We have 
$$\frac{d}{dt}\log Z_{dd}^{(\eps)}(\Phi_t)=\frac{d}{dt}\log\Qdet{K_t}=\frac12\frac{d}{dt}\log\det\tilde K_t,$$ and
\begin{eqnarray*}
\det(\tilde K_t+\delta \tilde S_t)&=&\det \tilde K_t\det(1+\delta \tilde S_t\tilde K_t^{-1})\\
&=& \det \tilde K_t\left(1+\delta\Tr(\tilde S_t\tilde K_t^{-1})+O(\delta^2)\right)\\
&=&\det \tilde K_t\left(
1+\delta\sum_{u_1,u_2} \tilde S_t(u_1,u_2)\tilde K_t^{-1}(u_2,u_1)+O(\delta^2)\right),
\end{eqnarray*}
so that $$\frac{d}{dt}\log\det \tilde K_t = \sum_{u_1,u_2} \tilde S_t(u_1,u_2)
\tilde K_t^{-1}(u_2,u_1).$$
Here $u_1,u_2$ run over all rows of $\tilde K_t^{-1}$.

Suppose without loss of generality that only $A_1$ changes.
Then $\tilde S_t(u_1,u_2)=0$ unless $u_1,u_2$ correspond to the
vertices of an edge of zipper $E_1$,
that is, $u_1u_2=w_{i}^{(p)}b_{i}^{(q)}$ or $u_1u_2=b_{i}^{(p)}w_{i}^{(q)}$
for an edge $w_{i}b_{i}$ of $E_1$, and $p,q\in\{1,2\}$. 
In this case $\tilde S_t(u_1,u_2)$ is equal to $\tilde K(u_1,u_2)$ (the edge weight) times either $A_1'(p,q)$ or $(A_1^{-1})'(p,q)$ depending on which direction
$u_1u_2$ crosses $E_1$. 

The sum of the two contributions from an edge $e_{i}=w_{i}b_{i}$ and its reverse is then 
\be\label{edgecontribution}
\tilde K_t(b_{i},w_{i})\sum_{p,q=1,2}\frac{dA_1(pq)}{dt}\tilde K_t^{-1}(b_{i}^{(q)},w_{i}^{(p)})+
\frac{dA_1^{-1}(pq)}{dt} \tilde K_t^{-1}(w_{i}^{(q)},b_{i}^{(p)}).
\ee

For simplicity (and without loss of generality, by path invariance) let us assume that the path $\gamma_1$
is polygonal with slope $\pm 1$, and for each segment of $\gamma_1$ the zipper $E_1$ consists of a zig-zag path, alternately horizontal and vertical in one of the directions
NE,NW,SW,SE (as in Figure \ref{Geps} for the SW path). Suppose moreover that each
zig-zag segment of $E_1$ has length which is an even number of lattice spacings (except perhaps the last segment),
that is, an even horizontal and even vertical length.

Let us then compute the contribution to the sum (\ref{edgecontribution}) for a given segment of $\gamma_1$.
Suppose first that it is oriented northeast.
Let $\{(w_i,b_i)\}_{i=1}^{4k}$ be the corresponding edges of $E_1$. 
We group these edges into packets of four consecutive edges, each packet consisting of two horizontal and two vertical edges.

Consider a horizontal edge $wb$ with white vertex of type $W_0$, in the same packet as a 
horizontal edge $w'b'$ with white vertex of type $W_1$.
By Lemma \ref{KDA} below, when $wb$ is not within $O(\eps)$ of the boundary,
$$K_t^{-1}(w,b)=A_1(\frac14 I+\eps\Re(F^\dag_+(z)+F_-(z,z))+O(\eps^2))$$
and 
$$K_t^{-1}(w',b')=A_1(\frac14 I+\eps\Re(F^\dag_+(z)-F_-(z,z))+O(\eps^2)).$$
Summing these two contributions gives 
$$\frac12A_1 + 2A_1\eps\Re F^\dag_+(z)+O(\eps^2).$$
For these edges in the reversed orientations the sum is the $q$-conjugate of this
(recall that $A_1^*=A_1^{-1}$ when $\det A_1=1$):
$$\frac12A_1^{-1} + 2\eps\Re F^{\dag*}_+(z)A_1^{-1}+O(\eps^2).$$
The contribution to the sum(\ref{edgecontribution}) is then
$$\frac12(A_1^{-1})'A_1+\frac12A_1'A_1^{-1} + 
\eps\left((A_1^{-1})'A_1\Re F^\dag_+(z)+A_1'\Re F^{\dag*}_+(z)A_1^{-1}\right)+O(\eps^2).$$
The leading term vanishes:
$$\frac12(A_1^{-1})'A_1+\frac12A_1'A_1^{-1}=0,$$ 
leaving the term of order $\eps$ and a negligible error.

For the two vertical edges $wb$ and $w'b'$ in the packet we have
similarly
$$K_t^{-1}(w,b)=A_1(\frac{i}4 I+i\eps\Im(F^\dag_+(z)+F_-(z,z))+O(\eps^2))$$
and 
$$K_t^{-1}(w',b')=A_1(\frac{i}4 I+i\eps\Im(F^\dag_+(z)-F_-(z,z))+O(\eps^2)).$$
Summing these two contributions gives 
$$\frac{i}2A_1 + A_1 i\eps\Im F^\dag_+(z)+O(\eps^2).$$
the net contribution to (\ref{edgecontribution}) from these two edges and their reverses is then 
(we must multiply by $-i$ which is the entry $\tilde K(w,b)=\tilde K(w',b')$, and the leading terms cancel as before)
$$\eps\left((A_1^{-1})'A_1\Im F^\dag_+(z)+A_1'\Im F^{\dag*}_+(z)A_1^{-1}\right)+ O(\eps^2).$$

Summing the contributions for the packet of four edges we get
$$(A_1^{-1})'A_1\Im(F^\dag_+(z)dz) + A_1'\Im(F^{\dag*}_+(z)dz)A_1^{-1})+O(\eps^2)$$
where we used the notation ``$dz$" to represent $\eps(1+i),$ the displacement from the beginning
of the packet to the end of the packet.

Now consider the case of a northwest segment of $\gamma_1$.
The only change is the contribution for the paired vertical edges $wb$ and $w'b'$. These 
are
$$K_t^{-1}(w,b)=A_1(-\frac{i}4 I+i\eps\Im(F^\dag_+(z)+F_-(z,z))+O(\eps^2))$$
and 
$$K_t^{-1}(w',b')=A_1(-\frac{i}4 I+i\eps\Im(F^\dag_+(z)-F_-(z,z))+O(\eps^2)).$$
Summing these and multiplying by $i$ the edge weight and the appropriate matrix we get
$$-\eps\left((A_1^{-1})'A_1\Im F^\dag_+(z)+A_1'\Im F^{\dag*}_+(z)A_1^{-1}\right)+O(\eps^2).$$
When added to the horizontal contribution, the net contribution for four edges is
$$(A_1^{-1})'A_1\Im(F^\dag_+(z)dz) + A_1'\Im(F^{\dag*}_+(z)dz)A_1^{-1})+O(\eps^2)$$
where ``$dz$" now represents the displacement $\eps(-1+i)$.

In a similar manner the other two possible directions of segments of $\gamma_1$ also
contribute 
$$(A_1^{-1})'A_1\Im(F^\dag_+(z)dz) + A_1'\Im(F^{\dag*}_+(z)dz)A_1^{-1})+O(\eps^2),$$
where ``$dz$" represents $\eps$ times $1+i,-1+i,-1-i,1-i$ according to the direction of the segment being
$NE,NW,SW,SE$ respectively. 

When $wb$ is within $O(\eps)$
of the boundary this formula
must be modified: the leading terms are no longer of modulus $1/4$. However as mentioned
in the comments after Lemma \ref{KDA} below,
we can choose a local trivialization of the bundle near $wb$
(by isotoping the zipper out of the way) to see that the leading terms
in $K^{-1}(w,b)$ and $K^{-1}(b,w)$
are replaced by the same constant $C_{bw}$ and thus still cancel as before. 
The subleading terms of order $\eps$ only differ from the above when $b$ is
within $o(1)$ of 
the boundary, and so when summed these boundary-error terms contribute a negligible amount.

In the limit $\eps\to0$ 
and the sum becomes the imaginary part of contour integral (the $O(\eps^2)$ term, when summed 
over the path, is at most $O(\eps)$ and drops out)
\be\label{int}
\frac{d}{dt}\log\det \tilde K_t = \Tr\left\{\frac{dA_1^{-1}}{dt}A_1\Im\left(\int_{\gamma_1}F^\dag_+(z)\, dz\right)+
\frac{dA_1}{dt}\Im\left(\int_{\gamma_1}F^{\dag*}_+(z)\, dz\right)A_1^{-1}\right\}.\ee
These are contour 
integrals of analytic functions depending only on the conformal type of the surface,
along the path $\gamma_1$. Thus the limit is conformally invariant.
\hfill$\square$

\section{$K^{-1}$ and discrete analyticity}

Let $U$ and $U_\eps$ be as in section \ref{setup}.
The goal of this section is to determine the asymptotic form of $K^{-1}(b,w)$
for adjacent vertices $b,w$ (Lemma \ref{KDA} below). In the case of trivial bundle this was worked out in
\cite{K.ci}. The proof there applies essentially without change to the case of a nontrivial flat bundle.
We give here an overview of the results of \cite{K.ci} and then indicate how the proofs change in the presence of a bundle.

\subsection{Discrete analytic functions}

Let $V_\eps, E_\eps, F_\eps$ be the vertices, edges and faces of $U_\eps$,
where $F_\eps$ includes the outer face $f_0$ and the face $f_i$ inside the $i$th
boundary component.
A {\bf discrete analytic function} \cite{Duffin} is a function $u+iv$, where 
$u:V_\eps\to\R$ and
$v:F_\eps\to \R$ 
satisfy the {\bf discrete Cauchy Riemann equations}:
for an edge $e=x_1x_2$, 
$$
u(x_2)-u(x_1) = v(f_1)-v(f_2),
$$
where $f_1,f_2$ are the two faces adjacent to edge $e$, and $f_1$ is the face to the
left when $e$ is traversed from $x_1$ to $x_2$.
This can be written succinctly as
\be\label{CR}
*du=dv
\ee
where $d$ represents the difference operator and $*$ is the ``rotation by $90^\circ$".

On non-boundary edges this is equivalent to  the ``discrete Cauchy-Riemann equations"
\begin{equation}\label{CReqns}\begin{array}{ccc}u_x&=&v_y\\
u_y&=&-v_x\end{array}
\end{equation}
on, respectively, horizontal and vertical non-boundary
edges, where the partial derivatives
represent discrete differences: $u_x((x,y)(x+\eps,y))=u(x+\eps,y)-u(x,y)$ and 
$u_y((x,y)(x,y+\eps))=u(x,y+\eps)-u(x,y).$ On boundary edges these
still hold if one interprets $v$ as being constant just outside each boundary
component. 

A function which satisfies (\ref{CReqns}) except at some subset of 
edges $e_1,\dots,e_k$ is 
said to be {\bf discrete meromorphic} with poles at $e_1,\dots,e_k$. In this case
the defect of the CR equations defines the {\bf residue}: the residue
for a horizontal edge $x_1x_2$ where $x_2=x_1+(\eps,0)$ and $f_1,f_2$ are the adjacent faces as in (\ref{CR}) is 
$u(x_2)-u(x_1)-v(f_1)+v(f_2)$ and the residue for a vertical edge $x_1x_2$ where $
x_2=x_1+(0,\eps)$ is
$i(u(x_2)-u(x_1)-v(f_1)-v(f_2))$. Note that the residue is either real (if the edge is horizontal)
or pure imaginary (if the edge is vertical). Also note that the residue is defined even for boundary edges.

If $u+iv$ is a discrete analytic function then both $u$ and $v$ are discrete
harmonic: 
\be\label{harm}
4u(p) = u(p+(\eps,0))+u(p+(0,\eps))+u(p-(\eps,0))+u(p-(0,\eps))
\ee
and likewise for $v$. This follows from summing (\ref{CR})
for the four edges coming out of a vertex (for $u$) or the four edges surrounding
a face (for $v$). If $u+iv$ is meromorphic with a pole at $e$ of residue $c\in\R$ then
$u$ is not harmonic at the vertices of $e$ and $v$ is not harmonic
at the faces adjacent to $e$. The Laplacian of $u$ at the vertices of $e$
is $\pm c$, depending on whether the vertex is the right or left endpoint of $e$. 
The Laplacian of $v$ is $\pm c$ at the upper, resp. lower face. Similar
equations hold for imaginary residues (at vertical edges).

If $U$ is a multiply-connected planar domain (or Riemann surface),
a {\bf discrete analytic section} of a flat bundle on $U_\eps$ is a section
which is locally a
discrete analytic function in any local trivialization of the bundle. 

\subsection{Kasteleyn matrix}

Recall the definition of the bipartite graph $\G_\eps$ from section \ref{setup}.
The black vertices of $\G_\eps$ are
$B=B_0\cup B_1$, where $B_0$ are vertices of $U_\eps$ and $B_1$ are
faces of $U_\eps$. The white vertices are $W=W_0\cup W_1$ where $W_0$ are horizontal edges
of $U_\eps$ and $W_1$ are vertical edges. 

Let $K_\eps$ be the Kasteleyn matrix for $\G_\eps$
whose rows index white vertices and columns index black vertices,
with $K_\eps(w,b)=1,i,-1,-i$ according to whether $b$ is adjacent and 
$E,N,W,$ or $S$ of $w$.

If we avoid the boundary, then $K_\eps$ acting on functions on $B$ is the 
discrete $\partial_{\bar z}$
operator in the sense that 
$u+iv$ is discrete analytic function on $U_\eps$ if and only if
$u+iv$, considered as a function on $B$ (that is, $u$ on $B_0$ and $iv$ on $B_1$)
is in the kernel of $K_\eps$. Taking the boundary
values into account, a function $u+iv$
in the kernel of $K_\eps$ must also satisfy $v=0$ on the large faces $f_i$ and the outer face, and satisfy $u=0$ at $x_0$, but can have poles at the $e_i$.

More generally, suppose that $(\G_\eps,U)$ is equipped with a flat bundle,
and $K_\eps$ the associated Kasteleyn matrix (whose entries
are obtained by multiplying the above
weights by the parallel transports $\phi_e$). Then, away from the boundary, a section is discrete analytic if and only if it is in the kernel of $K_\eps$. 

\begin{lemma}\label{Kinvmeromorphic}
As a function of $b$, $K_\eps^{-1}(b,w)$ is a discrete meromorphic (matrix-valued)
section with poles at the $e_j$ and at $w$, and zeros at $b=x_0$ and $b=f_j$ for all $j$. The pole at 
$w$ has residue $I$ or $iI$ according to $w\in W_0$ or $w\in W_1$. There is a unique section with these properties.
\end{lemma}

\begin{proof} Fix $w$; the equations $\sum_bK_\eps(w',b)K_\eps^{-1}(b,w)=\delta_{w,w'}I$ 
are linear equations for $K_\eps^{-1}(b,w)$, one for each $w'$.
At $w'\ne w$ they are the discrete CR equations. For $w'$ a boundary edge or 
$w'$ adjacent to $x_0$, they correspond to 
the CR equations if we extend $K_\eps^{-1}(b,w)$ to be zero at $b=f_i$ and $b=x_0$. 
The condition on the residue of the pole at $w$ is determined by $\sum_b
K_\eps(w,b)K_\eps^{-1}(b,w)=I$.
The uniqueness follows from invertibility of $K_\eps$.
\end{proof}
 
\subsection{Green's function}
The function $K_\eps^{-1}$ for the trivial line bundle on $\G_\eps$ is related to the Green's function $G$ 
of the standard Laplacian on $U_\eps$ and the Greens function $G^*$
on the dual graph $U_\eps^*$ as follows.
\begin{lemma}[\cite{K.ci}, Lemma 9]\label{KdiffG}
We have
$$K_\eps^{-1}(b,w) = \left\{\begin{array}{ll}G(w+\frac{\eps}2,b)-G(w-\frac{\eps}2,b)&w\in W_0, b\in B_0\\
-i(G^*(w+\frac{\eps}2 i,b)-G^*(w-\frac{\eps}2 i,b))&w\in W_0, b\in B_1\\
G^*(w+\frac{\eps}2,b)-G^*(w-\frac{\eps}2,b)&w\in W_1,b\in B_1\\
-i(G(w+\frac{\eps}2 i,b)-G(w-\frac{\eps}2 i,b))&w\in W_1,b\in B_0\end{array}\right.
$$
where $G,G^*$ is the Greens function for $U_\eps,U^*_\eps$ respectively with respectively
Neuman, Dirichlet boundary conditions.
\end{lemma}

With the appropriate definitions of discrete derivatives $\frac{\partial}{\partial w_x}$ and $\frac{\partial}{\partial w_y}$,
we can rewrite this as
$$K_\eps^{-1}(b,w) = \left\{\begin{array}{ll}\frac{\partial G(w,b)}{\partial w_x}&w\in W_0, b\in B_0\\
-i\frac{\partial G^*(w,b)}{\partial w_y}&w\in W_0, b\in B_1\\
\frac{\partial G^*(w,b)}{\partial w_x}&w\in W_1,b\in B_1\\
-i\frac{\partial G(w,b)}{\partial w_y}&w\in W_1,b\in B_0\end{array}\right..
$$

In \cite{K.ci} the asymptotics of $K^{-1}$ is written in terms of the continuous Green's function, as follows.
Let $\tilde g(u,v)$ be the analytic function of $v$ whose real part is the Dirichlet Green's function $g(u,v)$. 
Let $\tilde g^*(u,v)$ be the analytic function of $v$ whose real part is the Neumann Green's function $g^*(u,v)$.
Define\footnote{Note that this differs from the definition in \cite{K.ci}
by a factor of $4$; a factor of $2$ is conventional and one is due to the difference in choice of coordinates: here the lattice step for the dimer model
is $\eps/2$, not $\eps$.} $$F_+(u,v)=\frac{\partial\tilde g(u,v)}{\partial u}=
\frac12\left(\frac{\partial\tilde g(u,v)}{\partial u_x}-i\frac{\partial\tilde g(u,v)}{\partial u_y}\right)
$$
and 
$$F_-(u,v)=\frac{\partial\tilde g(u,v)}{\partial\bar u}=
\frac12\left(\frac{\partial\tilde g(u,v)}{\partial u_x}+i\frac{\partial\tilde g(u,v)}{\partial u_y}\right).$$
These are analytic functions of $u,v$ and $\bar u,v$ respectively.
It is not hard to show\footnote{For the upper half plane $\tilde g(u,v)=-\frac1{2\pi}\log\frac{u-v}{\bar u-v}$
and $\tilde g^*(u,v)=-\frac1{2\pi}\log(u-v)(\bar u-v).$ For a general simply connection domain $U$ with Riemann
map $\phi$ to $U$ we have $\tilde g(u,v)=-\frac1{2\pi}\log\frac{\phi(u)-\phi(v)}{\overline{\phi(u)}-\phi(v)}$ and
$\tilde g^*(u,v)=-\frac1{2\pi}\log(\phi(u)-\phi(v))(\overline{\phi(u)}-\phi(v))$.
}
that $$\frac{\partial\tilde g^*(u,v)}{\partial u}=F_+(u,v)\qquad \frac{\partial\tilde g^*(u,v)}{\partial\bar u}=- F_-(u,v),$$
that is, the difference is only a sign from $\tilde g$. 
Then for $w,b$ close to $u,v$ respectively with $u\ne v$ we have up to errors of order $O(\eps^2)$
$$K_\eps^{-1}(b,w) = \left\{\begin{array}{ll}
\eps\Re(F_+(u,v)+F_-(u,v))&w\in W_0, b\in B_0\\
\eps i\Im(F_+(u,v)+F_-(u,v))&w\in W_0, b\in B_1\\
\eps \Re(F_+(u,v)-F_-(u,v))&w\in W_1,b\in B_1\\
\eps i\Im(F_+(u,v)-F_-(u,v))&w\in W_1,b\in B_0\end{array}\right..
$$

Note that for fixed $w\in W_0$, $K^{-1}(w,b)$ is indeed analytic as a function of $v$, except at $v=u$.
Similarly for $w\in W_1$. 

In the case $w,b$ are within $O(\eps)$ of each other (for example adjacent), we define
$$F_+^\dag(u) = \lim_{v\to u} \left(F_+(u,v)-\frac1{2\pi(v-u)}\right).$$
Then 
\be\label{K4}
K_\eps^{-1}(b,w) = K^{-1}_{\eps,\Z^2}(b,w)+\left\{\begin{array}{ll}
\eps\Re(F^\dag_+(u)+F_-(u,u))&w\in W_0, b\in B_0\\
\eps i\Im(F^\dag_+(u)+F_-(u,u))&w\in W_0,b\in B_1\\
\eps\Re(F^\dag_+(u)-F_-(u,u))&w\in W_1, b\in B_1\\
\eps i\Im(F^\dag_+(u)-F_-(u,u))&w\in W_1,b\in B_0\end{array}\right..
\ee

Here note that $F^\dag_+(u)$ is analytic but $F_-(u,u)$ is not in general.

These formulas were shown to hold in the case of trivial bundle but the proof applies in the 
case of a flat $\SL$ connection as well; the relevant Green's functions are the
Green's function of the Laplacian 
for the connection $\Phi$, see \cite{K.bundles}; it is defined
by $$\Delta f(v) = \sum_{v'\sim v} f(v)-\phi_{v'v} f(v')$$ where the sum is over nearest
neighbors $v'$ of $v$.
The Green's functions $G$ and $G^*$ are
the inverse Laplacian operators on $U_\eps$ or $U^*_\eps$ respectively with the 
appropriate boundary conditions.
For generic $\Phi$ these Laplacians are invertible 
(see \cite[Theorem 9]{K.bundles} for an expression for their determinants) 
and so $G$ and $G^*$ are well-defined.

Since the Green's function $G_0$ on $\Z^2$ satisfies
$G_0((0,0),v)=-\frac14$ for any neighbor $v$ of $(0,0)$
and $G_0((0,0),(0,0))=0$ \cite{Spitzer},
by (\ref{K4}) we have the following result.

\begin{lemma}\label{KDA}
For an edge $wb$ with $w\in\W_0$, crossing the zipper $E_1$ at location $z\in U$ we have, up to $O(\eps)^2$ additive error,
$$K^{-1}(b,w) =\left\{\begin{array}{ll}A_1(\pm\frac14I + \eps\Re(F^\dag_+(z)+F_-(z,z)))&w\in W_0,~b=w\pm\eps/2\\
A_1(\mp\frac{i}4I + i\eps\Im(F^\dag_+(z)+F_-(z,z)))&w\in W_0,~b=w\pm i\eps/2\\
A_1(\pm \frac{1}4I +\eps \Re(F^\dag_+(z)-F_-(z,z)))&w\in W_1, b=w\pm\eps/2\\
A_1(\mp\frac{i}4I + i\eps\Im(F^\dag_+(z)-F_-(z,z)))&w\in W_1,~b=w\pm i\eps/2\\
\end{array}\right..
$$
\end{lemma}

When $z$ is near the boundary of $U$ (within $O(\eps)$) we have a weaker estimate.
We first isotope the zipper so
that it terminates at some other point not close to $z$; then the local behavior of $K^{-1}$ near $z$
is that $K^{-1}(b,w)=C_{bw}+O(\eps)$ for a constant $C_{bw}$ depending only on the distance
of edge $bw$ to the boundary and its orientation, see \cite{K.ci}. Now if we move the zipper back so that
it passes through $bw$ then $K^{-1}(b,w)=A_1 C_{bw}+O(\eps)$.

\section{Peripheral curves}\label{peripheral}

Theorem \ref{CI} applies to multiply connected domains with piecewise smooth boundaries.
If we are interested in simply-connected domains with punctures,
one can argue as follows.
A {\bf peripheral} curve is one which surrounds a single boundary component or puncture. 

On a multiply connected domain, the above proof shows that for any simple curve the 
probability that it occurs as a curve in the lamination of a random double-dimer cover is
conformally invariant. This applies to both peripheral and nonperipheral curves.
In the case when the boundary components shrink to points,
the number of peripheral curves surrounding them tends to infinity.
However the nonperipheral curves in a lamination are bounded (almost surely)
in number and
topological complexity. Thus one can make sense of the limiting probability of a nonperipheral curve.

The proof of Theorem \ref{CI} above shows that the computation of $F(\Phi)$ only depends on the
Green's function of the domain $U$. In particular if we introduce a finite number
of holes in $U$, with small diameter $<\delta$, which are not close to any existing
boundary or zipper then along any zipper the Green's function $G$ changes by an amount 
tending to zero with $\delta$. Since $F(\Phi)$ is an integral of $G$, 
$F(\Phi)$ also changes by an amount tending to zero with $\delta$. 
Thus the probabilities of homotopy classes of
double-dimer loops in $U$ are close to those in $U'$ 
(taking into account only the homotopy
classes relative to the boundary of $U$, not the new holes). 

That is, the double dimer loops are insensitive to the addition of small holes in the domain (except of course for those small loops which come close to or
intersect the removed holes). In particular if one wants to measure the probability of a homotopy class
of nonperipheral loop in the complement of a finite number of points, one can simply remove
small disks of radius $\delta$ around those points and compute the probability
in the limit as $\delta\to 0$.

\section{Other graphs}\label{othergraphs}

Theorem \ref{CI} applies to other graphs $U_\eps$ as well,
on condition that they conformally approximate $U$.

Let $U$ be as in section \ref{setup} and
let $U_\eps$ be a sequence of graphs which conformally approximate $U$
as $\eps\to0$, in the sense that the mesh size (diameter of the largest face)
tends to zero, and discrete harmonic functions on $U_\eps$ converge
to continuous harmonic functions on $U$. Equivalently,
the simple random walk on $U_\eps$ converges to the (time-rescaled) Brownian
motion on $U$. 

From such a graph $U_\eps$ one can define a bipartite
graph $G_\eps$ as above,
whose white vertices are the edges of $U_\eps$ and black vertices are the 
vertices and faces of $U_\eps$, see \cite{KPW}. The Kasteleyn matrix
is the adjacency matrix with rows indexing the white vertices and columns
indexing the black vertices, and (in general, complex-valued)
signs chosen so that around each face (which is a quadrilateral)  
the alternating product of signs is $-1$, that is,
if the face is a quadrilateral $ABCD$ with signs
$\sigma(AB),\dots,\sigma(DA)$), then
$\frac{\sigma(AB)\sigma(CD)}{\sigma(BC)\sigma(DA)}=-1$.

One natural way to choose the signs is to make them complex numbers
depending on the geometry of the embedding of $\G_\eps$, as follows.
For an edge $wb$ where $b$ is a vertex of $U$ let $K(w,b)=\frac{b-w}{|b-w|}$, that is, 
the unit modulus complex number in the direction of the edge $wb$. 
If $b$ is a face of $U$ let $K(w,b)$ be the unit complex number perpendicular to the
edge $w$ and in the direction from the edge to the face $b$.
This leads to the Kasteleyn weighting criterion that 
the four weights of a quadrilateral face $a,b,c,d$ satisfy $ac/bd=-1$.

The notion of discrete analytic function on $\G_\eps$ given by
(\ref{CR}) (but not (\ref{CReqns})) extends to this more general
setting, as does the harmonicity of the real and imaginary parts (equation (\ref{harm}),
with the usual combinatorial Laplacian on $U$).
Lemmas \ref{Kinvmeromorphic} and \ref{KdiffG} have simple analogs in this
more general setting:
$K^{-1}(b,w)$ is discrete meromorphic with a pole of residue $e^{i\theta}$
at $w$, when the edge is oriented in the direction $\theta$ (the residue is naturally
a $1$-form, so depends on a choice of orientation), and poles at any other removed edges of $U_\eps$. 
In terms of the Green's function, when $b$ is a vertex of $U$ we have
$$K^{-1}(b,w)=\frac{\partial G(w,b)}{\partial w}dw$$
where again $\frac{\partial}{\partial w}$ represents a difference operator on $U_\eps$.
A similar formula with the conjugate Green's function holds when $b$ is a face of $U$.

In Lemma \ref{KDA} we have the same formula except that the constants $1/4$ and $i/4$
are replaced by edge-dependent quantities (but they still cancel in (\ref{edgecontribution})). The subleading terms are the same by our hypothesis that
$U_\eps$ conformally approximates $U$.
We are led to exactly the same integral (\ref{int}) as before.

\section{Annulus}\label{annulus}

The case of a cylindrical annulus is particularly easy to compute because we can
explicitly diagonalize the relevant Kasteleyn matrices. 
Let $\G$ be the graph obtained from a rectangle $[0,2n]\times[1,m]\subset\Z^2$
by identifying for each $j$ vertex $(0,j)$ with vertex $(2n,j)$. This is a 
bipartite graph with $2nm$ vertices. 
We use a flat connection with monodromy $M$ around the 
circumference; thus a double-dimer configuration with $k$ noncontractible cycles will have weight
proportional to $(\Tr M)^k$. Since the weight only depends on $\Tr M$,
there is no loss of generality in taking $M$ to be a diagonal matrix
$M=\left(\begin{matrix}\lambda&0\\0&\lambda^{-1}\end{matrix}\right).$ Thus we need only consider 
a line bundle (see section \ref{diagonalmtx}).

Let $a$ satisfy $a^{2n}=\lambda$. We put parallel transport $a$ on horizontal edges $(x,y)(x+1,y)$ and $1$ on vertical edges.
To make a  Kasteleyn matrix $K$ multiply all vertical edge weights with $i=\sqrt{-1}$.
If we assume $n$ is odd then there are no additional signs needed.

The eigenvectors of $K$ are then of the following form.
Let $z,w$ satisfy $z^{2n}=1$ and $w^{2m+2}=1$. Then for vertices $(x,y)\in[0,2n]\times[1,m]$ the function
$$f_{z,w}(x,y) =z^x (w^y-w^{-y})$$ is an eigenvector with eigenvalue $$az+\frac1{az}+i(w+\frac1w).$$
Letting $z$ run over $2n$-th roots of unity
and $w$ run over the $2m+2$-th roots of $1$ with positive imaginary part 
(there are $m$ of these)
we have all $2nm$ independent eigenvectors. 

The determinant of $K$ is then, letting $w=e^{\pi ik/(m+1)}$ for $k=1,\dots,m$,
\begin{eqnarray*}
\det K &=& \prod_{z^{2n}=1}\prod_{k=1}^m\left(az+\frac1{az}+2i\cos\frac{\pi k}{m+1}\right)\\
&=&\prod_{z^{2n}=1}\prod_{k=1}^m\frac{(az-\alpha_k)(az-\beta_k)}{az}
\\
&=&\prod_{k=1}^m\frac{(\lambda-\alpha_k^{2n})(\lambda-\beta_k^{2n})}{\lambda}
\end{eqnarray*}
where $\alpha_k,\beta_k=i(-\cos\theta\pm\sqrt{1+\cos^2\theta})$ with 
$\theta=\frac{\pi k}{m+1}.$
Let $\alpha_k$ be the smaller (in modulus) of the two roots.
We can write this as
$$\det K=\left(\prod_k -\beta_k^{2n}\right)\left( \prod_k(1+\frac{|\alpha_k|^{2n}}{\lambda})(1+\lambda|\alpha_k|^{2n})\right)$$
since $\alpha_k\beta_k=1$ and $n$ is odd. Since $|\alpha_k|<1$, the terms in the second
product are negligible for large $n$ 
unless $|\alpha|\approx 1$, that is, except when $\theta\approx\pi/2$.
If $\theta=\frac{\pi}{2}+\eps$
for small $\eps$ then $|\alpha|^{2n}=e^{-2n|\eps|+O(n\eps^3)}.$
Suppose $m$ is even (the case $m$ odd is similar, see below). Take $k=\frac{m}2 + j$ so that $\theta=\frac{\pi k}{m+1} = \frac{\pi}2+\frac{\pi(2j-1)}{2(m+1)}$, and
$|\alpha_k|^{2n}=e^{-\frac{n\pi |2j-1|}{m+1}+O(nj/m^3)}.$
Taking $n,m$ large with $n/m=\tau$ fixed and  $q=e^{-\tau\pi}$ we have
\begin{eqnarray*}
\det K &\doteq& N\prod_{j\in\Z}(1+q^{|2j-1|}\lambda)(1+q^{|2j-1|}\lambda^{-1})\\
&=&N\prod_{\substack{j=1\\j\text{ odd}}}^\infty(1+q^{2j}+q^{j}X)^2,
\end{eqnarray*}
where $X=\lambda+\lambda^{-1}=\Tr(M)$
and where $N$ is a normalizing factor independent of $\lambda$.

\begin{figure}[htbp]
\includegraphics[height=2in]{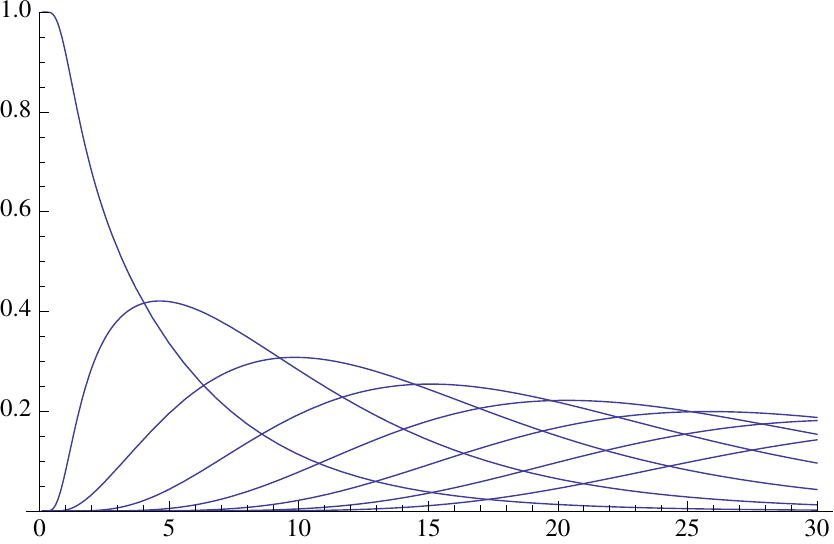}
\caption{\label{loopprobs} Probabilities of $0$ to $7$ loops as a function of $1/\tau=m/n,$
the ratio of height to circumference of the cylinder, in the case $m$ even.}
\end{figure}

This gives the probability generating function for the number of loops to be (when $m$ is even)
$$\sum_{k=0}^\infty \Pr(k \text{ loops})X^k  = \prod_{\substack{j=1\\j\text{ odd}}}^\infty
\frac{(1+q^{j}X+q^{2j})^2}{(1+q^{j}+q^{2j})^2}.$$
See Figure \ref{loopprobs}.

Similarly one finds, in the case $m$ is odd, the probability generating function to be 
$$\sum_{k=0}^\infty \Pr(k \text{ loops})X^k  = \frac{2+X}{3}\prod_{\substack{j=2\\j\text{ even}}}^\infty
\frac{(1+q^{j}X+q^{2j})^2}{(1+q^{j}+q^{2j})^2}.$$

\section{Loops surrounding two points}\label{twofaces}
\begin{figure}[htbp]
\includegraphics[height=2.in]{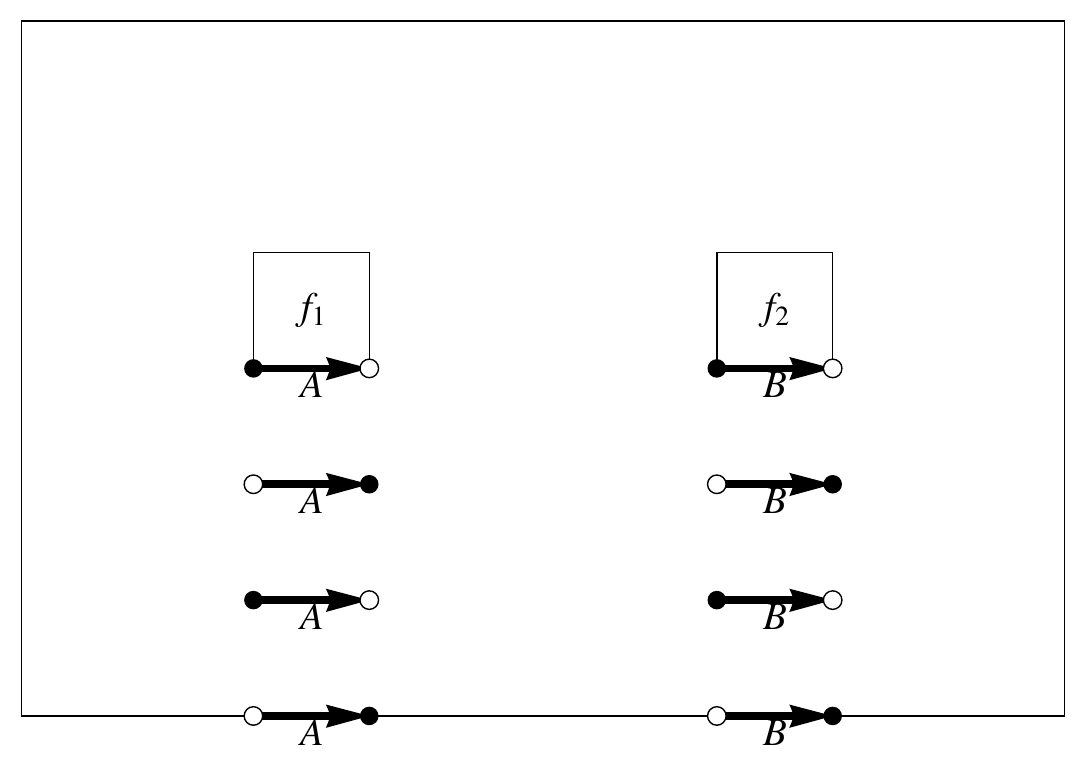}
\caption{\label{twofacefig}}
\end{figure}

We show here how to compute the distribution for the number of loops
surrounding two faces in the upper half plane grid $\G=\Z^2\cap\{y>0\}$.
Let $f_1,f_2$ be two faces, and take a flat $\SL$-bundle with
monodromy $A\in\SL$ on a zipper from the boundary to $f_1$ and $B\in\SL$ on
a zipper from the boundary to $f_2$ as in Figure \ref{twofacefig}.

We let $A=\left(\begin{matrix}1&\eps\\0&1\end{matrix}\right)$ and
$B=\left(\begin{matrix}1&0\\\eps&1\end{matrix}\right)$ for small $\eps$.
Then $\Tr(AB)=2+\eps^2$.

We then have 
$$Z=\Qdet K=\sum_{k=0}^\infty C_k(1+\frac{\eps^2}{2})^k,$$
where $C_k$ is the weighted sum 
(weighted by $2^c$, $c$ being the total number of loops)
of configurations with $k$ loops
surrounding both $f_1$ and $f_2$. (Note that loops surrounding only one of $f_1$ or $f_2$ are counted correctly since $\Tr A=2=\Tr B$.)

The coefficient of $\eps^2$ in the expansion around $\eps=0$ of $Z$
is then $\frac12$ of the expected number of loops.
Higher coefficients give higher moments. 

Let us compute the coefficient of $\eps^2$. We have $K=\left(\begin{matrix}
0&M\\M^*&0\end{matrix}\right),$ where $M$ is the matrix with $\SL$ entries and with
rows
indexing white vertices and columns indexing black vertices.
Let $M_0$ be the corresponding matrix for the trivial connection;
the double $\tilde M_0$ is a direct sum of two copies of $K_0$, the standard scalar
Kasteleyn matrix for the upper half plane. 
Now $\Qdet K=\sqrt{\det\tilde K} =\det\tilde M$ and so 
$$Z=\det \tilde M = \det(\tilde M_0+\eps S)=
\det \tilde M_0+\eps^2\sum_{w_1b_1,w_2b_2}(\tilde M_0)_{w_1^{1}w_2^{2}}^{b_1^{2}b_2^{1}}+O(\eps^4).$$
Here $w_1b_1$ is an edge of zipper $A$ and $w_2b_2$ is an edge of zipper $b$.
(The terms of first order in $\eps$ are zero, as are the terms of order
$\eps^2$ with both edges from the same zipper.)

However $\tilde M_0$ is just two copies of $K_0$, so  
$$(\tilde M_0)_{w_1^{1}w_2^{2}}^{b_1^{2}b_2^{1}}=(\det\tilde M_0)K_0^{-1}(w_1^1b_2^1)K_0^{-1}(w_2^2b_1^2).$$

In the scaling limit for the upper half plane, let $w_1b_1$ be a horizontal
edge near a point with complex coordinate $z_1$ and $w_2b_2$ be a horizontal edge near a point with complex coordinate $z_2$. 
Taking into account the orientations (Right or Left) 
of $w_1b_1$ and $w_2b_2$ we have
$K_0^{-1}(w_1^1b_2^1)K_0^{-1}(w_2^2b_1^2)=$
$$=\left\{\begin{array}{ll}
\frac{\eps^2}{\pi^2}\Re(\frac1{z_2-z_1}+\frac1{z_2-\bar z_1})\Re(\frac1{z_1-z_2}+\frac1{z_1-\bar z_2})&
\text{if }(w_1^1b_2^1,w_2^2b_1^2)\text{ is RR}\\
-\frac{\eps^2}{\pi^2}\Im(\frac1{z_2-z_1}+\frac1{z_2-\bar z_1})\Im(\frac1{z_1-z_2}-\frac1{z_1-\bar z_2})&\text{if }(w_1^1b_2^1,w_2^2b_1^2)\text{ is RL}\\
-\frac{\eps^2}{\pi^2}\Im(\frac1{z_2-z_1}-\frac1{z_2-\bar z_1})\Im(\frac1{z_1-z_2}+\frac1{z_1-\bar z_2})&\text{if }(w_1^1b_2^1,w_2^2b_1^2)\text{ is LR}\\
\frac{\eps^2}{\pi^2}\Re(\frac1{z_2-z_1}-\frac1{z_2-\bar z_1})\Re(\frac1{z_1-z_2}-\frac1{z_1-\bar z_2})&\text{if }(w_1^1b_2^1,w_2^2b_1^2)\text{ is LL}
\end{array}\right.
$$

Summing over the four possibilities, the sum is a Riemann sum for the integral
$$\frac{Z}{\det\tilde M_0}[\eps^2] = -\frac{2}{\pi^2}\int_{\gamma_1}\int_{\gamma_2}
\Re\left(\frac{1}{(z_1-z_2)^2}+\frac{1}{(z_1-\bar z_2)^2}\right)\,dy_1\, dy_2$$
$$=-\frac{2}{\pi^2}\log\left|\frac{z_1-z_2}{z_1-\bar z_2}\right|$$
up to errors going to zero with $\eps$.

Thus we have
\begin{theorem}For the double-dimer model on $\eps\Z^2\cap\{y>0\}$,
the expected number of 
loops surrounding both of the points $z_1,z_2$ converges as $\eps\to0$ to 
$$-\frac{4}{\pi^2}\log\left|\frac{z_1-z_2}{z_1-\bar z_2}\right|.$$
\end{theorem}

This result can also be obtained using the results of \cite{K.ci}
and fact that the double-dimer loops are the contours of the
height function difference of two independent uniform dimer covers.

An analogous computation for a chordal path can be made:
take $b,w$ two vertices on the boundary of the upper half plane grid 
$\G=\eps\Z^2\cap\{y\ge0\}$ with $b<w$. Add to $\G$ an edge $e$ connecting $b$ to $w$.
Take a double dimer cover of $\G$, and condition on the event that edge $e$ is part of a nontrivial loop.
This loop will be a
``chordal" double-dimer path from $b$ to $w$ in the upper half plane. 

\begin{theorem}
In the scaling limit on the upper half plane, the probability that the chordal
double-dimer path from $b$ to $w$ with $b<w\in\R$ passes left of a point $z$ is the harmonic
function of $z$ with boundary values $1$ between $b$ and $w$ and zero 
elsewhere.
\end{theorem}

\end{document}